\documentclass[12pt]{article}
\usepackage{amssymb}
\usepackage{graphicx}
\usepackage{amsfonts}
\usepackage{latexsym}
\usepackage{amsthm}
\usepackage{amsmath}

\usepackage{amssymb, amscd}

\usepackage{url}

\newtheorem{problem}{Problem}
\newtheorem{theo}[problem]{Theorem}
\newtheorem{rem}{Remark}
\newtheorem{prob}[problem]{Problem}
\newtheorem{defin}[problem]{Definition}
\newtheorem{prop}[problem]{Proposition}
\newtheorem{cor}[problem]{Corollary}
\newtheorem{lema}[problem]{Lemma}
\newtheorem{exam}[problem]{Example}

\newtheorem{conj}[problem]{Conjecture}

\begin{document}
\date{July 12, 2013}
 \title{Illumination complexes, $\Delta$-zonotopes, \\ and the polyhedral curtain theorem}

\author{{Rade  T.\ \v Zivaljevi\' c\footnote{Supported by Grants 174020 and 17434 of the
Serbian Ministry of Education and Science.}}\\ {\small Mathematical Institute}\\[-2mm] {\small SASA, Belgrade}
\\[-2mm]{\small rade$@$mi.sanu.ac.rs}}

\maketitle
\begin{abstract}
{\em Illumination complexes} are examples of `flat polyhedral
complexes' which arise if several copies of a convex polyhedron
(convex body) $Q$ are glued together along some of their common
faces (closed convex subsets of their boundaries). A particularly
nice example arises if $Q$ is a $\Delta$-zonotope (generalized
rhombic dodecahedron), known also as the dual of the difference
body $\Delta - \Delta$ of a simplex $\Delta$, or the dual of the
convex hull of the root system $A_n$. We demonstrate that the
illumination complexes and their relatives can be used as
`configuration spaces', leading to new `fair division theorems'.
Among the central new results is the `polyhedral curtain theorem'
(Theorem~\ref{thm:prva}) which is a relative of both the `ham
sandwich theorem' and the `splitting necklaces theorem'.
\end{abstract}

\section{Introduction}\label{sec:introduction}

`Polyhedral curtain theorem' is a combinatorial relative of both
the classical `ham sandwich theorem' and Alon's `splitting
necklaces theorem'. In general, a `polyhedral curtain' or a
`polyhedral wall' in $\mathbb{R}^d$ is a polyhedral set in
$\mathbb{R}^d$ of dimension $(d-1)$ which is homeomorphic to
$\mathbb{R}^{d-1}$ and which separates $\mathbb{R}^d$ into two
connected components (polyhedral half-spaces). In this paper we
deal mainly with conical polyhedral curtains generated by
polyhedral spheres in $\mathbb{R}^d$.

\begin{defin}\label{def:pol-surface}
A conical, polyhedral hypersurface $D = {\rm cone}(a,\Sigma)$,
where $\Sigma\subset \mathbb{R}^d$ is a $(d-2)$-dimensional
polyhedral sphere ($PL$-sphere) and $a\notin \Sigma$, is called a
{\em polyhedral curtain} in $\mathbb{R}^d$.
\end{defin}

\begin{figure}[hbt]
\centering
\includegraphics[scale=0.40]{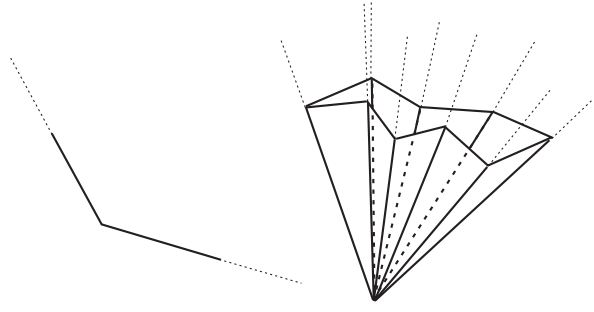}
\caption{Conical polyhedral curtains in $\mathbb{R}^2$ and
$\mathbb{R}^3$.} \label{fig:pol-curtain-2-3-dim}
\end{figure}

The simplest example of a polyhedral curtain is the union of two
rays in $\mathbb{R}^2$ emanating from the same point
(Figure~\ref{fig:pol-curtain-2-3-dim}).
Our main examples of polyhedral curtains in $\mathbb{R}^d$ are
$\Delta$-curtains in the sense of the following definition.
\begin{defin}
Let $\Delta = \mbox{\em\rm conv}\{a_0,a_1,\ldots, a_d\}\subset
\mathbb{R}^d$ be a non-degenerate simplex with the barycenter at
the origin. For each pair $\theta=(F_1, F_2)$ of complementary
faces of $\Delta$ there is a join decomposition $\Delta = F_1\ast
F_2$. Assuming that both $F_1$ and $F_2$ are non-empty let
$S^{d-2}_\theta = \partial(F_1)\ast
\partial(F_2)\subset \partial(\Delta)$ be an associated $(d-2)$-dimensional, polyhedral
sphere. A polyhedral hypersurface $H\subset \mathbb{R}^d$ is
called a $\Delta$-{\em curtain} if for some $\theta = (F_1,F_2)$
and $x\in \mathbb{R}^d$
\begin{equation}
H = x + C_\theta = x + \mbox{\rm cone}(S^{d-2}_\theta).
\end{equation}
\end{defin}
\begin{figure}[hbt]
 \hspace{3.5cm}
\includegraphics[scale=0.70]{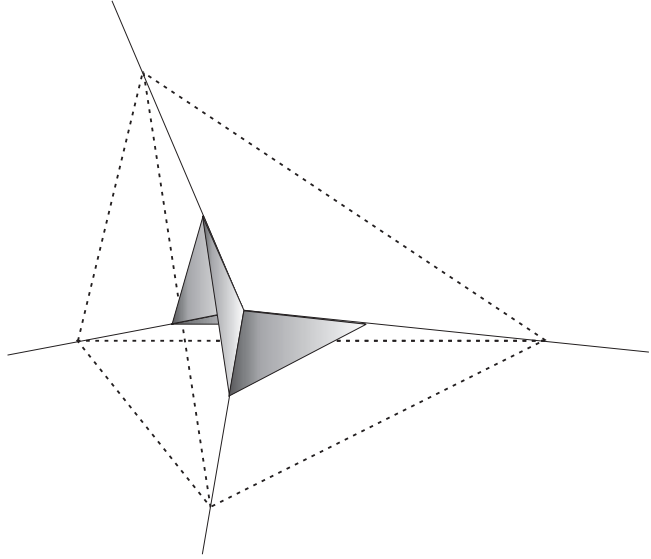} \caption{A fragment of
a $\Delta$-curtain in $\mathbb{R}^3$.} \label{fig:3d-curtain}
\end{figure}
The classical ham sandwich theorem claims that a collection of $d$
measurable sets in $\mathbb{R}^d$ admits a hyperplane bisector,
i.e.\ a hyperplane that simultaneously cuts them in halves of
equal measure. The following `polyhedral curtain theorem' is
formally a statement of similar nature where the role of
hyperplanes is played by $\Delta$-curtains.

\begin{theo}\label{thm:prva} (Polyhedral curtain theorem)
Suppose that $\Delta\subset \mathbb{R}^d$ is a simplex with the
barycenter at the origin. Let $\mu_1, \mu_2,\ldots, \mu_d$ be a
collection of continuous mass distributions (measures) on
$\mathbb{R}^d$. Then there exists a $\Delta$-curtain $H = H_{(x,
\theta)}$ which divides the space $\mathbb{R}^d$ into two
`half-spaces' $H^+$ and $H^-$ such that for each $j\in\{1,\ldots,
d\}$,
\[
\mu_j(H^+) = \mu_j(H^-).
\]
\end{theo}

\begin{rem}{\rm  In contrast with the ham sandwich theorem, Theorem~\ref{thm:prva} has a more combinatorial
flavor, in this respect it resembles Radon's and Tverberg theorem
\cite{Z04}. It will become clear from the proof that it is really
an offspring of the multidimensional splitting necklace theorem
\cite{Long-Ziv}, a higher dimensional generalization of the
celebrated splitting necklace theorem of Alon, \cite{Alo87,
alon-constructive}.
 }
\end{rem}

\section{Two dimensional case of Theorem~\ref{thm:prva}}
\label{sec:two-dim-case}

\begin{figure}[hbt]
\centering
\includegraphics[scale=0.40]{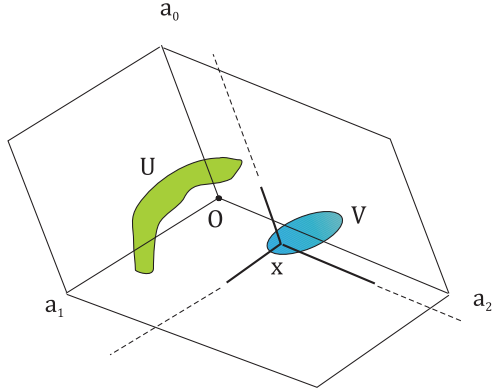}
\caption{A division with a moving fan.}
\label{fig:pol-curtain-hexagon}
\end{figure}

As a motivation for introducing $\Delta$-zonotopes
(Section~\ref{sec:delta-zonotopi}) and {\em flat (polyhedral)
complexes} (Section~\ref{sec:flat}), we outline the proof of the
two dimensional case of the {\em polyhedral curtain theorem}
(Theorem~\ref{thm:prva}), emphasizing the main ideas and
constructions. We will demonstrate that two (bounded) measurable
sets $U$ and $V$ in the plane (say the simple shapes depicted in
Figure~\ref{fig:pol-curtain-hexagon}) admit a fair division by a
planar polyhedral curtain, with the directions of the rays
prescribed in advance by a triangle $\Delta$. Without a loss of
generality we assume that $U\cup V \subset \Delta$.

\begin{figure}[hbt]
\centering
\includegraphics[scale=0.40]{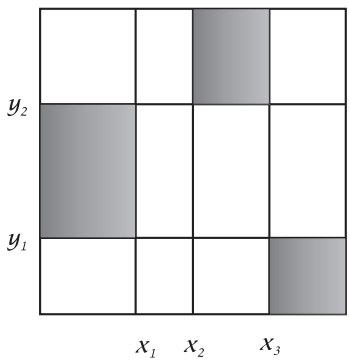}
\caption{Division of the square with axis-aligned cuts.}
\label{fig:kvadrat-podela}
\end{figure}

\medskip
Theorem~\ref{thm:prva} is a `fair division theorem', so like in
other results of this type one is supposed to identify the
associated `configuration space' of allowed (admissible)
divisions. For comparison, let us briefly review the construction
of the configuration space for the two dimensional splitting
necklace problem, \cite[Section~2]{Long-Ziv} (with $3$ vertical
and $2$ horizontal axis aligned cuts).  An admissible division of
the square $I^2$, corresponding to the chosen sequences $0=x_0\leq
x_1\leq x_2\leq x_3\leq x_4=1; 0= y_0\leq y_1\leq y_2\leq y_3 =
1$, is depicted in Figure~\ref{fig:kvadrat-podela}. It immediately
follows that the configuration space of all divisions of $I^2$ (of
the type $(3,2)$) is the product $\Delta^3\times \Delta^2$ of
simplices.

Assuming that there are two persons involved in the division, we
observe that there are $2^{12}$ possible scenarios for
distributing these rectangles (one possibility is depicted in
Figure~\ref{fig:kvadrat-podela}). As a consequence, the
configuration space $\Omega(3,2)$ of all divisions of the square
of the type $(3,2)$, together with all possible allocations of
rectangular pieces to two parties involved, is the union of
$2^{12}$ polyhedral cells (copies of $\Delta^3\times \Delta^2$).
These cells are glued together, along their boundaries, whenever
some of the elementary rectangles degenerate (for example if $x_i
= x_{i+1}$ or $y_j = y_{j+1}$), see \cite[Section~2]{Long-Ziv} for
more details.

In the two dimensional case of the polyhedral curtain theorem,
instead of the square, the basic shape (convex body) is a
centrally symmetric hexagon, associated to the triangle $\Delta =
{\rm conv}\{a_0,a_1,a_2\}$, Figure~\ref{fig:pol-curtain-hexagon}.
The axis-aligned cuts (of the square) are replaced by the
dissection of the hexagon $Q$ into three regions, determined by
the three rays emanating from the same point (the apex of the
associated fan). The rays are always assumed to be translates of
the basic system of rays (basic fan), generated by vectors $\{a_0,
a_1, a_2\}$ with the apex at the origin $O$. The apex of the
translated fan can be any point in the hexagon, including the
boundary, so the configuration space of all such division is the
hexagon $Q$ itself. There are $2^3=8$ possibilities to allocate
the three regions, obtained by this division, to the two
interested parties. As a result there are eight hexagons, each
associated to one of possible eight scenarios for the division.
These hexagons are glued together whenever one (or two) of the
three regions degenerates, this happens if the apex $x$ is on the
boundary of $Q$.

Summarizing, we observe that the configuration space encoding all
admissible divisions, together with the associated division
scenarios, is a cell complex obtained by gluing together eight
identical copies of a convex body (hexagon $Q$) along their
boundaries, following a specific gluing scheme. We denote this
configuration space by $\mathbb{A}(Q,\mathcal{F},S)$ (see also
Section~\ref{sec:flat}) where $\mathcal{F}$ is the associated fan
and $S = [2] =\{1,2\}$ the set of associated `colors'
(representing the parties involved in the division).

A moment's reflection shows that $\mathbb{A}(Q,\mathcal{F},[2])$
is homeomorphic to the $2$-sphere. Moreover, the involution on the
set $[2]$ (corresponding the parties interchanging their roles),
defines an involution on the set of cells and on the configuration
space $\mathbb{A}(Q,\mathcal{F},[2])$, which is easily identified
as the usual antipodal action.

\medskip
In order to complete the proof of the planar case of the
polyhedral curtain theorem we construct a `test map' $\phi :
\mathbb{A}(Q,\mathcal{F},[2]) \rightarrow \mathbb{R}^2$, testing
the fairness of the chosen division-allocation. More explicitly,
let $(x,s)\in \mathbb{A}(Q,\mathcal{F},[2])$, where $x$ is the
apex of the associated fan $\mathcal{F}=\{F_0, F_1, F_2\}$ and $s
: \mathcal{F} \rightarrow \{0,1\}$ is the associated allocation
function. Let $X_0 = \cup\{F_i \mid s(F_i) = 0\}$ and $X_1 =
\cup\{F_i \mid s(F_i) = 1\}$ be the associated polyhedral
`half-planes', the planar regions separated by the chosen
polyhedral curtain. Then by definition,
\[
\phi((x,s)) = (m(X_0\cap U) - m(X_1\cap U), m(X_0\cap V) -
m(X_1\cap V)).
\]
This function is easily shown to be continuous. By construction
$\phi(-z) = -\phi(z)$. Since the zeros of $\phi$ correspond to
fair divisions, the planar case of the polyhedral curtain theorem
follows as a consequence of the Borsuk-Ulam theorem.

\begin{exam}{\rm (communicated by Sini\v sa Vre\' cica)
It is interesting that one cannot in general guarantee the
existence of a bisecting $\Delta$-curtain for three measurable
sets $U, V, W$, even if the rotations (isometries) of the
$\Delta$-curtains are allowed. Indeed, the example is provided by
small discs centered at the vertices of the triangle $\Delta$. }
\end{exam}

\section{$\Delta$-zonotopes $R_\Delta$}
\label{sec:delta-zonotopi}

Generalizing and developing the ideas from
Section~\ref{sec:two-dim-case} we introduce $\Delta$-zonotopes as
higher dimensional analogues of the centrally symmetric hexagon,
used in the proof of the planar case of Theorem~\ref{thm:prva}.

\begin{defin}\label{def:delta-polytope}
Let $\Delta = {\rm conv}\{a_0, a_1, \ldots , a_d\}$ be a
non-degenerate simplex in $\mathbb{R}^d$ such that $a_0+\ldots +
a_d = 0$. The convex polytope, defined as the Minkowski sum,
\begin{equation}
 R_\Delta = [0,a_0] + [0, a_1] + \ldots + [0, a_d]
\end{equation}
is referred to as $\Delta$-zonotope. The simplex $\Delta$ is
called the {\em generating simplex} of the polytope $R_\Delta$.
Each polytope $R_\Delta$ has a `standard cubulation' $R_\Delta =
C_0\cup\ldots \cup C_d$ where $C_i$ is the parallelotope
($d$-parallelepiped) spanned by vectors $\{a_j\}_{j\neq i}$.
\end{defin}
\begin{figure}[hbt]
\centering
\includegraphics[width=6 cm]{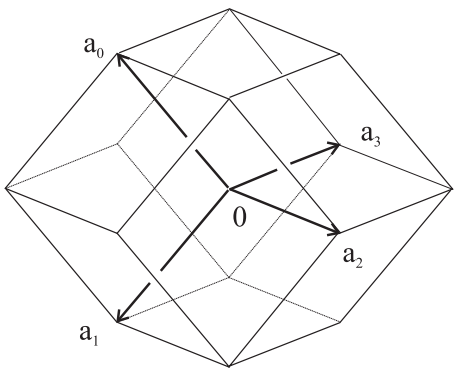}
\caption{Rhombic dodecahedron and its generating simplex.}
\label{fig:R12-3}
\end{figure}

Two dimensional $\Delta$-zonotopes are the regular hexagon and its
affine images. Figure~\ref{fig:R12-3} shows the rhombic
dodecahedron, the $3$-dimensional $\Delta$-zonotope, together with
its generating intervals (generating simplex). Recall that
according to the classic classification of the mathematician,
crystallographer, and mineralogist Evgraf Fedorov, rhombic
dodecahedron is one of the five space filling polytopes known as
{\em parallelohedra}. Rhombic dodecahedron is also the Voronoi
cell of the {\em face centered cubic lattice} in the $3$-space.

\medskip
Rhombic dodecahedron (and $\Delta$-zonotopes as its higher
dimensional analogues) have appeared in a very interesting problem
of Makeev \cite{Mak1} about universal covers of convex bodies of
diameter $1$. Recall that the related result in dimension $2$  is
one of the classics of the combinatorial geometry in the plane
(\cite{Pal}).

\begin{conj}\label{conj:Makeev}{\rm (V.V.~Makeev \cite{Mak1}) Every
convex body $K$ in $\mathbb{R}^d$ of diameter $\leq 1$ can be
covered by a translate of the convex body $R_\Delta$ for some
non-degenerate simplex $\Delta$ of diameter $\leq 1$.}
\end{conj}

The conjecture was subsequently established in dimension $d=3$,
independently by Makeev \cite{Mak2} (with mild assumptions on
$K$), by Hausel, Makai, and Sz\"{u}cs \cite{Ha-Ma-Sz}, and by
Kuperberg \cite{Kuper}. The original conjecture was formulated for
the dual $D_d = (\Delta - \Delta)^\circ$ of the difference body of
a simplex, rather than for the $\Delta$-zonotope $R_\Delta$. The
added importance of the difference body $\Delta - \Delta = {\rm
conv}\{e_i - e_j \mid 0 \leq n\}$ of the $d$-simplex $\Delta =
{\rm conv}\{e_i\}_{i=0}^n\subset \mathbb{R}^{n+1}$ stems from the
fact that it can be described as the convex hull of the root
system $A_{n}$. The following proposition shows the equivalence of
these definitions of generalized rhombic dodecahedra.

\begin{figure}[hbt]
\centering
\includegraphics[scale=0.50]{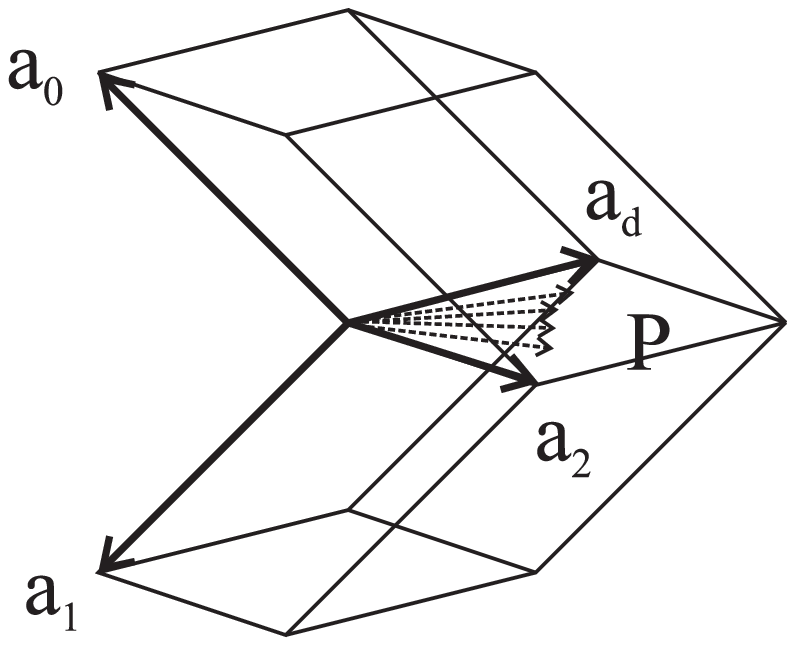}
\caption{$\Delta$-zonotope $R_\Delta$ is the dual of
$\Delta-\Delta$.} \label{fig:R12-general}
\end{figure}

\begin{prop}\label{prop:dual-difference-body}
Let $\Delta = {\rm conv}\{a_0, a_1, \ldots , a_d\}$ be a regular
simplex in $\mathbb{R}^d$ such that $a_0+\ldots + a_d = 0$. Let
$H_{i,j} = \{x\in \mathbb{R}^d \mid \vert \langle x, a_i - a_j
\rangle \vert \leq  \vert \langle a_i, a_i-a_j \rangle \vert \}$
be the region between two hyperplanes orthogonal to the edge $a_i
- a_j$, passing respectively through the vertices $a_i$ and $a_j$.
Then for some constant $\lambda >0$,
\begin{equation}\label{eqn:pos-constant}
\lambda (\Delta - \Delta)^\circ = \bigcap_{i\neq j} H_{i,j} =
R_\Delta .
\end{equation}
\end{prop}

\medskip\noindent
{\bf Proof:} (outline). Two parallelotopes $C_0$ and $C_1$ from
the standard cubulation of $R_\Delta$
(Definition~\ref{def:delta-polytope}) are depicted in
Figure~\ref{fig:R12-general}. Let $P$ be the parallelotope
generated by vectors $a_2,\ldots, a_d$ and let $H$ be the
associated hyperplane. Then $H$ is orthogonal to $a_0-a_1$ and the
supporting hyperplanes $a_0+ H$ and $a_1+H$ of $R_\Delta$ are
precisely the boundary hyperplanes of $H_{0,1}$. \hfill $\square$

\subsection{$\Delta$-zonotopes and quadrangulations}
\label{sec:quadrangulation}

$\Delta$-zonotopes appear, at least implicitly, in toric topology
in the context of standard quadrangulations (cubulations) of
simple polytopes, see \cite[Chapter 4]{Buch-Pan}. Here we
summarize one of this constructions, referring the reader to
\cite{Buch-Pan} for more details and other related information.
\begin{figure}[hbt]
\centering
\includegraphics[scale=0.60]{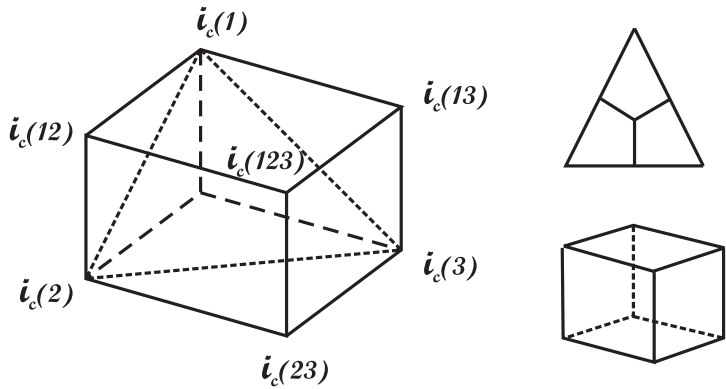}
\caption{The front complex of a cube.} \label{fig:front-and-back}
\end{figure}

Each non-empty subset $I = \{j_1,\ldots, j_k\}\subset [n]$ can be
associated both a vertex $a_I = (1/k)(e_{j_1}+\ldots +e_{j_k})$ of
the barycentric subdivision of the simplex $\Delta =
\Delta^{n-1}={\rm conv}\{e_1,\ldots, e_n\}$ and the vertex $b_I =
e_{j_1}+\ldots +e_{j_k}$ of the standard cube $I^n\subset
\mathbb{R}^n$. The correspondence $a_I \mapsto b_I$ is extended to
a piecewise linear embedding $i_c : \Delta^{n-1} \rightarrow I^n$
which turns out to be an isomorphism of $\Delta^{n-1}$ with the
`front face' $F_f(I^n)$ of the cube $I^n$
(Figure~\ref{fig:front-and-back}).

The front face $F_f(I^n)$ is formally described as the union of
all facets of $I^n$ which contain the vertex $b_{[n]}=e_1+\ldots
+e_n$. A moment reflection shows that there exists a natural
piecewise linear isomorphism $\nu : F_f(I^n) \rightarrow
R_{\Delta^{n-1}}$ of the front face $F_f(I^n)$ and the
$\Delta$-zonotope associated to $\Delta^{n-1}$.  Indeed, the
standard cubulation of $R_{\Delta^{n-1}}$
(Definition~\ref{def:delta-polytope}) is naturally brought into
one-to-one correspondence with the natural cubulation of the front
face of $I^n$.

Some of the properties of the piecewise linear isomorphisms $i_c$
and $\nu$, as well as their composition $I_c = \nu\circ i_c$,
\begin{equation}\label{eqn:isomorphism-I}
I_c : \Delta^{n-1} \stackrel{i_c}{\longrightarrow} F_f(I^n)
\stackrel{\nu}{\longrightarrow} R_{\Delta^{n-1}},
\end{equation}
are summarized in the following proposition.

\begin{prop}\label{prop:summary-I}
The map $I_c$ described by (\ref{eqn:isomorphism-I}) is a
piecewise linear isomorphism which preserves the associated
cubical decompositions (Figure~\ref{fig:front-and-back}). The
standard cubulation $R_{\Delta^{n-1}} = C_1\cup\ldots \cup C_n$
(Definition~\ref{def:delta-polytope}) induces a decomposition
\[
\partial (R_{\Delta^{n-1}}) = F_f(C_1)\cup\ldots \cup F_f(C_n)
\]
of the boundary of the $\Delta$-zonotope $R_{\Delta^{n-1}}$ into
the associated front faces. This decomposition is well-behaved
with respect to the isomorphism $I_c$ in the sense that
$I_c(F_i(\Delta^{n-1})) = F_f(C_i)$, where $\partial
(\Delta^{n-1}) = \cup_{i=1}^n~F_i(\Delta^{n-1})$ is the facet
decomposition of the simplex $\Delta^{n-1}$.
\end{prop}

\section{Flat polyhedral complexes}\label{sec:flat}

A {\em flat polyhedral complex} $E$ arises if several copies of a
convex polyhedron (convex body) $B$ are glued together along some
of their common faces (closed convex subsets of their boundaries).
By construction there is a `folding map' $p : E\rightarrow B$
which resembles the moment map from toric geometry. Examples of
flat polyhedral complexes include `small covers' and other locally
standard $\mathbb{Z}_2$-toric manifolds, but the idea can be also
traced back (at least) to A.D.~Alexandrov's `flattened convex
surfaces'.

A class of flat polyhedral complexes, modelled on a product of two
or more simplices $B = \Delta^{p_1}\times\ldots\times
\Delta^{p_k}$, was used in \cite{Long-Ziv} for a proof of a
multidimensional generalization of Alon's `splitting necklace
theorem' \cite{Alo87}. Developing this idea we demonstrated in
Section~\ref{sec:introduction} that some new classes of flat
polyhedral complexes naturally appear as `configuration spaces',
leading to new `fair division theorems'.

Here we introduce two classes of flat complexes, based on convex
bodies (polyhedra) in $\mathbb{R}^d$ and develop their theory as
far as it is needed for our central applications, including the
proof of Theorem~\ref{thm:prva}.

\subsection{Convex fans and illumination systems}
\label{sec:illumination-systems}

A convex fan $\mathcal{F}$ (\cite{Zig}) is tacitly assumed to have
the apex at the origin. Moreover, it is complete in the sense that
it covers the whole of $\mathbb{R}^d$. For a given convex fan
$\mathcal{F}$ in $\mathbb{R}^d$ let $\mathcal{F}_{max} = \{V_1,
V_2, \ldots, V_m\}$ be the associated collection of maximal cones
in $\mathcal{F}$. Since $\mathcal{F}$ is determined by
$\mathcal{F}_{max}$, we shall (with a mild abuse of the language)
often neglect the difference between $\mathcal{F}$ and
$\mathcal{F}_{max}$ and sometimes write simply $\mathcal{F} =
\{V_1, V_2, \ldots, V_m\}$. A translated fan is defined as $x +
\mathcal{F} = \{x+V_i\}_{i=1}^m $.

\medskip
There is a well known class of problems in classical combinatorial
geometry known as {\em illumination} or {\em visibility} problems,
see \cite{Mar-Solt}. Motivated by this, a convex cone $V\subset
\mathbb{R}^d$, with the apex at $a\in \mathbb{R}^d$ may be
interpreted as the region in $\mathbb{R}^d$ {\em illuminated} by a
light source (reflector) positioned at the point $a$. More
generally, a fan $\mathcal{F}_{max} = \{V_1, V_2, \ldots, V_m\}$
can be interpreted as an `illumination system' of essentially
non-overlapping light sources, associated with the cones $V_i$.

Given a finite set of labels or colors $S$, the {\em status} or
the {\em state} of the illumination system $\mathcal{F}_{max} =
\{V_1, V_2, \ldots, V_m\}$ is a function $f : [m] \rightarrow S$
describing the state of each of the individual light sources. For
example for $i\in [m]$, the associated value $f(i)$ may be
interpreted as the color of the light chosen for the cone $V_i$
from the given set $S$ of colors. In the simplest case when $S =
\{0, 1\}$, the two possibilities naturally correspond to the
status of the light source being turned on or turned off.

Summarizing, we observe that an illumination system is a pair
$(\mathcal{F}, S)$ where $\mathcal{F}$ is a complete fan of convex
cones while $S$ the set of admissible colors for each of the
individual light sources.

\subsection{Illumination complexes $\mathbb{A}(K, \mathcal{F}, S)$}
\label{sec:illumination}

Let us start with the following data:

\begin{itemize}
\item $K\subset \mathbb{R}^d$ is a convex polyhedron (or more
generally a convex body);

\item $\mathcal{F}$ is finite complex of convex cones in
$\mathbb{R}^d$ (a convex fan);

\item $S$ is finite set of `colors'.

\end{itemize}

The illumination complex $\mathbb{A}(K, \mathcal{F}, S)$ is so
designed to encode all possible ways of illuminating the body $K$
from some point $x\in K$ by the (translated) illumination system
$(\mathcal{F}, S)$. More precisely, a typical element of
$\mathbb{A}(K, \mathcal{F}, S)$ is a  pair $(x, f)$ where $x\in K$
and $f : [m] \rightarrow S$ is a status of the illumination system
$(\mathcal{F}, S)$.
\begin{figure}[hbt]
\centering
\includegraphics[width=11 cm]{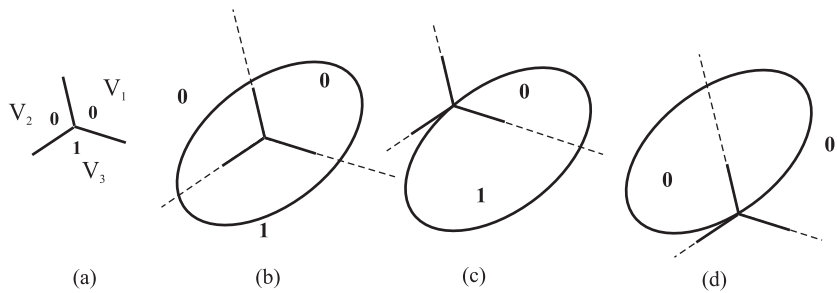}
\caption{Elements of the illumination complex $\mathbb{A}(K,
\mathcal{F}, S)$.} \label{fig:oval-2.eps}
\end{figure}
The choice of $x\in K$ determines where in $K$ the illumination
system should be positioned (translated) while the function $f :
[m]\rightarrow S$ describes the chosen inner status of each of the
individual light sources of the illumination system $(\mathcal{F},
S)$.

In the example depicted in Figure~\ref{fig:oval-2.eps}~(a), the
status function $f : [3] \rightarrow \{0,1\}$ dictates that the
light is turned off in cones $V_1$ and $V_2$ and turned on in
$V_3$. The three chosen positions of the point $x\in K$
(Figure~\ref{fig:oval-2.eps}~(b), (c) and (d)) illustrate typical
situations that may occur. In (b) all three light sources
illuminate a part of the interior of $K$. In (c) ${\rm
int}(V_2)\cap {\rm int}(K) = \emptyset $ which means that the
inner status of the light source associated with the cone $V_2$ is
no longer important and can be either $0$ or $1$. This observation
is a motivation for making the identification of points $(x,f)$
and $(x,g)$ where $g : [3] \rightarrow \{0,1\}$ is the status
function such that $f(i) = g(i)$ for $i\neq 2$.

After the informal description of the configuration space
$\mathbb{A}(K, \mathcal{F}, S)$ we are finally ready for a formal
definition.

\begin{defin}\label{def:A-kompleks}
The configuration space (cell complex) $\mathbb{A}(K, \mathcal{F},
S)$ is by definition the identification space $(K \times X)/
\sim$, where $X = Fun([m], S)$ is the set of all functions from
the index set $[m]$ (of $\mathcal{F}$) to $S$ and $(x,f)\sim
(y,g)$ if and only if, $x = y$ and $f(i) = g(i)$ for each $i\in
[m]$ such that ${\rm int}(V_i)\cap {\rm int}(K)\neq \emptyset$.
\end{defin}

\subsection{Modified illumination complex $\mathbb{B}(K, S)$}
\label{sec:modified-illumination}

The construction of the modified illumination complex
$\mathbb{B}(K, S)$ follows the same general pattern used in the
case of $\mathbb{A}(K, \mathcal{F}, S)$, with some important
differences. First of all $K$ is assumed to be a convex polytope
(in the previous case $K$ was allowed to be a convex body).
Secondly, for each $x\in K$ there is an associated fan
$\mathcal{F}_x = \mathcal{F}_x(K)$ (with the apex at $x$) where
$\mathcal{F}_x = \{{\rm cone}(x, F) \mid F \mbox{ \rm is a proper
face of } K\}$. Observe that if $x\in {\rm bd}(K)$, the dimension
of the cone ${\rm cone}(x,F)$ is $d$ if and only if $F$ is a facet
of $K$ such that $x\notin F$.

Summarizing, instead of a fan $\mathcal{F}$ which is prescribed in
advance, and the associated translated fan $x + \mathcal{F}$, here
we have a fan $\mathcal{F}_x$ with the moving apex $x$, generated
by the facets of $K$.

If $F(K) = \{F_i\}_{i=1}^m$ is the set of all facets of $K$ then
the status function describing the colors selected for each of the
associated light sources is a function $f : [m] \rightarrow S$.

\begin{defin}\label{def:B-kompleks}
The configuration space (cell complex) $\mathbb{B}(K, S)$ is by
definition the identification space $(K \times X)/ \sim$, where $X
= Fun([m], S)$ is the set of all functions from the index set
$[m]$ (indexing the facets of $\mathcal{F}$) to $S$ and $(x,f)\sim
(y,g)$ if and only if, $x = y$ and $f(i) = g(i)$ for each $i\in
[m]$ such that $x\notin F_i$.
\end{defin}

\subsection{The complex $\mathbb{B}(\Delta, S)$}
\label{sec:B-is-join}

\begin{theo}\label{thm:B-simplex}
Suppose that $\Delta = {\rm conv}\{a_0, a_1, \ldots, a_d\}\subset
\mathbb{R}^d$ is a non-degenerate simplex with the barycenter at
the origin, and let $S$ be a finite, non-empty set (of colors).
Then the associated $\mathbb{B}$-complex is homeomorphic to the
join of $(d+1)$ copies of the $0$-dimensional complex $S$,
$$
\mathbb{B}(\Delta, S) \cong S\ast S \ast \ldots \ast S = S^{\ast
(d+1)}
$$
\end{theo}

\begin{cor}\label{cor:B-simplex-2-boje}
If $S = \{0,1\}$ then the associated $\mathbb{B}$-complex
$\mathbb{B}(\Delta, [2])$ is isomorphic to the $d$-sphere, the
boundary $\partial \diamondsuit^d$ of the $d$-dimensional
cross-polytope $\diamondsuit^d$.
\end{cor}

\medskip\noindent
{\bf Proof of Theorem~\ref{thm:B-simplex}:}  Suppose that the
simplex $\Delta = {\rm conv}\{a_i\}_{i=0}^d$ has volume $1$,
$m(\Delta)=1$. Let us define a map,
\begin{equation}\label{eqn:iso}
\Phi : \mathbb{B}(\Delta, S) \rightarrow S^{\ast (d+1)}
\end{equation}
and show that it is an isomorphism of simplicial complexes. Let
$(x,f)\in  \mathbb{B}(\Delta, S)$. Let $x = \alpha_0 a_0 + \ldots
+ \alpha_d a_d$ be the barycentric decomposition $x$ and let $F_i$
be the facet of $\Delta$, opposite to the vertex $a_i$, $F_i =
{\rm conv}\{a_0,\ldots, \hat{a}_i, \ldots, a_d\}$.

\medskip\noindent
{\bf Lemma:} The volume of the `pyramid' $D_i(x):={\rm
conv}(\{x\}\cup F_i)$, with the apex at $x$ and the face $F_i$ as
the base, is equal to the barycentric coordinate $\alpha_i$. In
other words the knowledge of the volumes of all convex sets
$D_i(x)$ allows us to determine uniquely the position of the point
$x\in \Delta$.

\medskip
By definition $\Phi(x,f) = \alpha_0 f(0) + \alpha_1 f(1) + \ldots
+ \alpha_d f(d)\in S^{\ast(d+1)}$. It is not difficult to see that
this map is well defined and that it provides the desired
isomorphism between $\mathbb{B}(\Delta, S)$ and $S^{\ast (d+1)}$.
\hfill $\square$

\subsection{Well illuminated complexes}
\label{sec:well-illuminated}

Contrary to the case of the $\mathbb{B}$-complex
(Theorem~\ref{thm:B-simplex}), the corresponding
$\mathbb{A}$-complex $\mathbb{A}(\Delta, \mathcal{F}, S)$,
associated to the simplex $\Delta = {\rm conv}\{a_j\}_{j=0}^d$, is
quite irregular even if $\mathcal{F}$ is the fan $\mathcal{F}_0 =
\{{\rm cone\{F_i\}}\}_{i=0}^d$. The following definition clarifies
which properties of $\mathbb{A}$-complexes should be considered
`regular', at least from the point of view of intended
applications in this paper.

\begin{defin}\label{def:A-well-illminated}{\rm
Given a convex fan $\mathcal{F}$ in $\mathbb{R}^d$ and a finite
set $S$ of colors, we say that a convex body $K\subset
\mathbb{R}^d$ is $(\mathcal{F}, S)$-{\em well illuminated} if the
complex $\mathbb{A}(K, \mathcal{F}, S)$ is $(d-1)$-connected.  }
\end{defin}

\begin{defin}\label{def:fans} {\rm
If $Q$ is a convex polytope such that $0\in {\rm int}(Q)$ there
are two `tautological' fans associated with $Q$. The first is the
face fan,
\begin{equation}\label{eqn:face-fam}
\mathcal{F}_Q = \{{\rm cone}(F)\mid F \mbox{ {\rm is a proper face
of} } Q\},
\end{equation}
and the second is the normal fan $\nu(Q)$. }
\end{defin}

\begin{prob}\label{prob:A-tela}{\em
For a given pair $(\mathcal{F}, S)$ characterize convex bodies $K$
which are $(\mathcal{F}, S)$-well illuminated in the sense of
Definition~\ref{def:A-well-illminated}. As a first step it would
be interesting to know examples of well illuminated bodies,
especially if $\mathcal{F}$ is one of the fans $\mathcal{F}_Q$ and
$\nu(Q)$, associated to a convex body $Q$, and $S=\{0,1\}$ is a
$2$-element set.  }
\end{prob}

\subsection{The complex $\mathbb{A}(R_\Delta, \mathcal{F}_\Delta, S)$}
\label{sec:R-polytopes}

In this section we show that the centrally symmetric hexagon (in
the plane), the {\em rhombic dodecahedron} (in the $3$-space) and
their higher dimensional analogues $R_\Delta$ are all
$(\mathcal{F}_\Delta, S)$-well illuminated, if $\Delta$ is the
associated `generating simplex'
(Definition~\ref{def:delta-polytope}).

\begin{theo}\label{thm:randevu}
Let $\Delta = {\rm conv}\{a_0, a_1, \ldots , a_d\}$ be a
non-degenerate simplex in $\mathbb{R}^d$ with the barycenter at
the origin and let $R_\Delta$ be the associated $\Delta$-polytope
(Definition~\ref{def:delta-polytope}). Let $S$ be a finite,
non-empty set of colors. Let $\mathcal{F}_\Delta = \{{\rm
cone}(F)\mid F \mbox{ {\rm is a proper face of} } \Delta\}$ be the
fan generated by the facets of the simplex $\Delta$. Then,
\begin{equation}\label{eqn:randevu}
\mathbb{A}(R_\Delta, \mathcal{F}_\Delta, S) \cong S\ast S \ast
\ldots \ast S = S^{\ast (d+1)}.
\end{equation}
Moreover, this isomorphism is equivariant with respect to any
permutation of the colors $\pi : S\rightarrow S$. In particular if
$S=\{0,1\}$, the corresponding $\mathbb{A}$-complex is
homeomorphic to the $d$-sphere, the boundary of the
cross-polytope, $\mathbb{A}(R_\Delta, \mathcal{F}_\Delta, [2])
\cong
\partial(\diamondsuit^d)$.
\end{theo}

\medskip\noindent
{\bf Proof: }   Let $(x,f)\in \mathbb{A}(R_\Delta,
\mathcal{F}_\Delta, S)$. Define $F_i = {\rm conv}\{a_0, \ldots,
\hat{a}_i, \ldots, a_d\}$ as the facet of $\Delta$, opposite to
the vertex $a_i$, and let $\{{\rm cone}(F_i)\}_{i=0}^d$ be the
collection of maximal cones in $\mathcal{F}_\Delta$. The required
isomorphism,
\begin{equation}\label{eqn:rama}
\Psi : \mathbb{A}(R_\Delta, \mathcal{F}_\Delta, S) \rightarrow
S\ast S \ast \ldots \ast S = S^{\ast (d+1)}.
\end{equation}
is defined by,
\begin{equation}\label{eqn:isomorphism-psi}
\Psi(x,f) = (1/\alpha)[\alpha_0 f(0) + \alpha_1 f(1) + \ldots +
\alpha_d f(d)]
\end{equation}
where $\alpha_i = \alpha_i(x,f) = m((x+{\rm cone}(F_i))\cap
R_\Delta)$ is the volume of the region in $R_\Delta$ illuminated
by the cone $x+{\rm cone}(F_i)$ and $\alpha = \alpha_0+\ldots
+\alpha_d = m(R_\Delta)$.

\begin{lema}\label{lema:equivariance}{\rm
The map $\Psi$ is equivariant with respect to the group of all
permutations of the set $S$ of colors. }
\end{lema}

\medskip\noindent
{\bf Proof of Lemma~\ref{lema:equivariance}:} This is a
consequence of the fact that the cell $C_f\subset
\mathbb{A}(R_\Delta, \mathcal{F}_\Delta, S)$, associated to a
function $f : \{0, \ldots, d\} \rightarrow S$, is mapped to the
corresponding cell $D_f = \{f(0)\}\ast\ldots \ast \{f(d)\}\subset
S^{\ast (d+1)}$. In other words $\Psi$ is the equivariant
extension of the map $\psi : R_\Delta \rightarrow \Delta^d$,
defined by
\begin{equation}\label{eqn:basic-map}
\psi(x) = \alpha_0a_0+\ldots +\alpha_d a_d,
\end{equation}
where $\alpha_i$ have the same meaning as in the equation
(\ref{eqn:isomorphism-psi}). \hfill $\square$

\medskip\noindent

In the following lemma we show that the map $\Psi$ defined by
(\ref{eqn:isomorphism-psi}) is a monomorphism.

\begin{lema}\label{lema:one-to-one}{\rm
If $(x,f), (y, g)\in \mathbb{A}(R_\Delta, \mathcal{F}_\Delta, S)$
are distinct points then $\Psi(x,f)\neq \Psi(y,g)$. }
\end{lema}

\medskip\noindent
{\bf Proof of Lemma~\ref{lema:one-to-one}:} We can assume that
$f=g$, indeed if $\alpha_i(x,f)\neq 0 \neq \alpha_i(y,g)$ and
$f(i)\neq g(i)$ then clearly $\Psi(x,f)\neq \Psi(y,g)$. Let us
assume that $x$ is in the interior of $R_\Delta$. Since
$\mathcal{F}_\Delta$ is a complete fan, $y-x\in {\rm cone}(F_i)$
for some $i$. As a consequence,
\begin{equation}\label{eqn:relation-containment}
(y+{\rm cone}(F_i))\cap {\rm int}(R_\Delta)) \varsubsetneq (x+{\rm
cone}(F_i))\cap {\rm int}(R_\Delta)
\end{equation}
which implies that $\alpha_i(y,f)< \alpha_i(x,f)$. A slight
extension of this argument applies also in the case when both
points $x$ and $y$ are on the boundary of $R_\Delta$. Indeed, a
least one of the cones $V_i := x+{\rm cone}(F_i)$ such that
$y-x\in V_i$ has the property $V_i\cap {\rm
int}(R_\Delta))\neq\emptyset$ and the relation
(\ref{eqn:relation-containment}) is again satisfied. \hfill
$\square$

\medskip\noindent

\begin{figure}[hbt]
\centering
\includegraphics[scale=0.60]{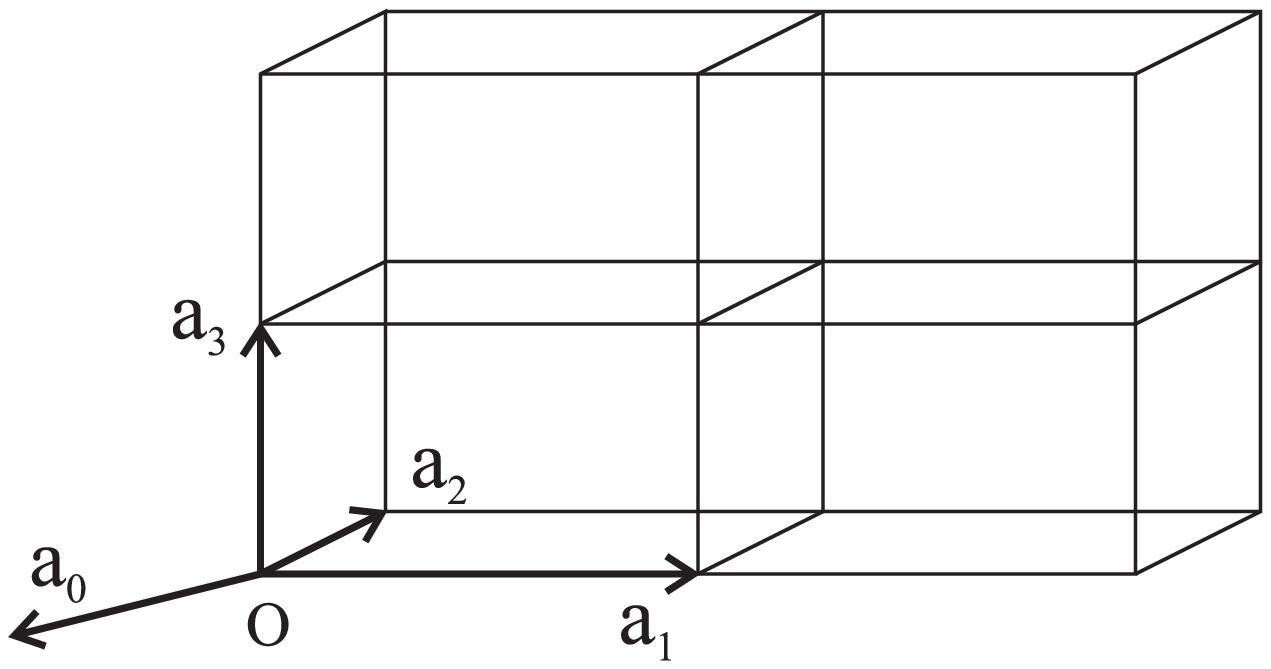}
\caption{} \label{fig:cigle-1}
\end{figure}

\noindent Let us now show that $\Psi$ is an epimorphism. Since
${\rm Image}(\Psi)\subset S^{\ast (d+1)}$ is a closed set, it is
sufficient to show that the image is everywhere dense in $S^{\ast
(d+1)}$.

\begin{lema}\label{lema:onto}{\rm
If $z\in S^{\ast (d+1)}$ has the representation $z = \beta_0 z_0 +
\beta_1 z_1 + \ldots + \beta_d z_d \in S^{\ast (d+1)}$ such that
$\beta_i\neq 0$ for each $i$, then for some $(x,f)\in
\mathbb{A}(R_\Delta, \mathcal{F}_\Delta, S), \, \alpha_i(x,f) =
\beta_i$. }
\end{lema}

\medskip\noindent
{\bf Proof of Lemma~\ref{lema:onto}:} It is sufficient to show
that the map $\psi : R_\Delta \rightarrow \Delta^d$, defined by
(\ref{eqn:basic-map}) is an epimorphism. Let $\phi = \psi\circ I_c
: \Delta^d \rightarrow \Delta^d$ where $I_c : \Delta_d \rightarrow
R_\Delta$ is the map defined by (\ref{eqn:isomorphism-I})
(Proposition~\ref{prop:summary-I}). By inspection of the
Figure~\ref{fig:cigle-1} we observe that $\alpha_i(x)=0$ if $x\in
F_f(C_i)$, the front face of $C_i$, where $R_\Delta =
C_0\cup\ldots\cup C_d$ is the standard cubulation of $R_\Delta$.
From here and Proposition~\ref{prop:summary-I} it immediately
follows that the map $\phi' : \partial(\Delta^d) \rightarrow
\partial(\Delta^d)$ is a homotopic to the identity map. By the
standard argument, used in the proof of Brouwer's fixed point
theorem, it follows that ${\rm int}(\Delta^d)\subset {\rm
Image}(\phi)$ which in turn implies that $\psi$ is also an
epimorphism. \hfill $\square$

\begin{cor}\label{cor:A-B-comparison}{\rm In light of Theorem~\ref{thm:B-simplex}
the complexes $\mathbb{A}(R_\Delta, \mathcal{F}_\Delta, S)$ and
$\mathbb{B}(\Delta, S)$ are equivariantly homeomorphic. Moreover,
this homeomorphism is naturally defined in terms of the associated
illumination systems.}
\end{cor}

\section{Fair divisions by polyhedral cones}

The general {\em configuration space/test map}-scheme
(CS/TM-scheme) was used long before it was recognized in
\cite{Ziv-Guide}, as a key organizing principle for applying the
topological methods on problems of combinatorial nature. According
to this scheme \cite{Ziv-Guide, Z04}, the problem can be
classified by the nature of its {\em configuration space} (the
space of all reasonable `candidates for the solution'), and the
nature of topological principles involved in its solution. The
emphasis on `configuration spaces' led to the systematization of
tools for their analysis and construction and, as a useful
byproduct, seemingly distant problems were recognized as neighbors
or (topological) `genetic relatives', see \cite{Matousek,
Ziv-Guide, Z04} for details and examples.

Following the CS/TM-scheme, the illumination complexes
$\mathbb{A}(K, \mathcal{F}, S)$ and $\mathbb{B}(K, S)$ were
designed as {\em configuration spaces} suitable for applications
to `envy-free', `fair division' or `consensus division' theorems,
where two or more parties are involved in dividing an object
following the rules prescribed in advance (see \cite{Long} for an
introduction and first examples, and \cite{Z04} for an overview).
In this section we formulate and prove these fair division
results.

As in Section~\ref{sec:flat} the simplex $\Delta = {\rm
conv}\{a_j\}_{j=0}^d\subset \mathbb{R}^d$ is assumed to be
non-degenerate with barycenter at the origin, and
$\mathcal{F}_\Delta = \{{\rm cone}(F_i)\}_{i=0}^d$ is the
associated `face-fan' determined by the facets $F_i$ of $\Delta$
(Definition~\ref{def:fans}).

\begin{theo}\label{thm:fair-prva} Choose positive integers $n,q$  such
that $q = p^k$ is a power of a prime $p\geq 2$ and let $d =
n(q-1)$. Let $\mu_1,\mu_2,\ldots, \mu_n$ be a collection of
continuous measures on $\mathbb{R}^d$ with finite support. Then
there exists a translate  $x + \mathcal{F}_\Delta = \{x + {\rm
cone}(F_i)\}_{i=0}^d$ of the face fan of $\Delta$ and an
allocation function $\theta : \{0,1,\ldots, d\}\rightarrow
\{1,\ldots, q\}$ such that for each $\mu_i$ and all
$j=1,2,\ldots,q$,
\begin{equation}\label{eqn:fair-prva}
\sum_{\theta(\nu)=j} \mu_i(x + {\rm cone}(F_\nu)) =
(1/q)\mu_i(\mathbb{R}^d) .
\end{equation}
\end{theo}

\medskip\noindent
{\bf Proof:} Without loss of generality (by a homothetic
enlargement of the simplex $\Delta$) we can assume that the
supports of all measures $\mu_i$ are contained in $\Delta$. Let
$R_\Delta$ be the $\Delta$-zonotope associated to $\Delta$
(Section~\ref{sec:delta-zonotopi}) and let $\mathbb{A}(R_\Delta,
\mathcal{F}_\Delta, S)$ be the associated illumination complex
(Section~\ref{sec:illumination}) where $S = [q]$ is the chosen set
of `colors'.

Recall that an element $(x,f)\in \mathbb{A}(R_\Delta,
\mathcal{F}_\Delta, S)$ is a pair of a point $x\in R_\Delta$ and a
status function $f : \{x+V_i\}_{i=0}^d \rightarrow [q]$ of the
illumination system $x + \mathcal{F} = \{x + V_i \mid 0\leq i\leq
d\}$, where $V_i = {\rm cone}(F_i)$.

For each pair of indices $i\in [n]$ and $j\in [q]$, define
$\alpha_{i,j} = \alpha_{i, j}(x,f)$ by
\begin{equation}\label{eqn:apha-i-j}
\alpha_{i,j} = \alpha_{i, j}(x,f) = \sum_{\theta(\nu)=j} \mu_i(x +
{\rm cone}(F_\nu)).
\end{equation}
The function $\alpha_{i,j}: \mathbb{A}(R_\Delta,
\mathcal{F}_\Delta, S) \rightarrow \mathbb{R}$ is well-defined and
continuous. The associated matrix valued function
$$
\Phi : \mathbb{A}(R_\Delta, \mathcal{F}_\Delta, S) \rightarrow
{\rm Mat}_{n\times q}(\mathbb{R}), \, \Phi(x,f) = A =
(\alpha_{i,j})
$$
is $\Sigma_q$-equivariant with respect to the action which
permutes the columns of the matrix $A = (\alpha_{i,j})$. Observe
that the unique fixed point of this action is the matrix $O =
(\omega_{i,j})$ such that $\omega_{i,j} =
(1/q)\mu_i(\mathbb{R}^d)$ for each $i$ and $j$. Observe that ${\rm
Image}(\Phi)\subset L$ where $L\subset {\rm Mat}_{n\times
q}(\mathbb{R})$ is a $\Sigma_q$-invariant affine subspace define
by equations,
$$
\sum_{j=1}^q~x_{i,j} = \mu_i(\mathbb{R}^d)\quad \mbox{ {\rm for
each} }\, i=1,\ldots, n.
$$
The dimension of this space is $n(q-1)$ and if
(\ref{eqn:fair-prva}) is never satisfied there arises a
$\Sigma_q$-equivariant map,
\[
\Psi : \mathbb{A}(R_\Delta, \mathcal{F}_\Delta, S) \rightarrow
S(L)
\]
where $S(L)$ is the $[n(q-1)-1]$-dimensional unit sphere in $L$.
Since by Theorem~\ref{thm:randevu} the complex
$\mathbb{A}(R_\Delta, \mathcal{F}_\Delta, S)$ is
$(d-1)$-connected, the existence of such a map would be in
contradiction with the following general Borsuk-Ulam type theorem
(Theorem~\ref{thm:Volovikov}).

\begin{theo}\label{thm:Volovikov}
{\rm (\cite{Vol, Sar, Long-clanak, Ziv-Guide, Matousek, Long})}
Suppose that $p\geq 2$ is a prime number and let $G =
\mathbb{Z}_p^{\times k}$ be the elementary abelian group of order
$p^k$. Assume that $X$ and $Y$ are fixed-point free $G$-spaces
such that $X$ is $n$-connected and $Y$ is a $n$-dimensional
sphere. Then there does not exist a $G$-equivariant map $f : X
\rightarrow Y$.
\end{theo}

Theorem~\ref{thm:fair-prva} relied in essential way on the
properties of the complex $\mathbb{A}(R_\Delta,
\mathcal{F}_\Delta, S)$. A $\mathbb{B}(K,S)$-counterpart of that
result is the following theorem.

\begin{theo}\label{thm:fair-druga} Choose positive integers $n,q$  such
that $q = p^k$ is a power of a prime $p\geq 2$ and let $d =
n(q-1)$. Let $\Delta\subset \mathbb{R}^d$ be a non-degenerate
simplex. Assume that $\mu_1,\mu_2,\ldots, \mu_n$ are continuous
measures on $\mathbb{R}^d$ with finite support ${\rm
Supp}(\mu_i)\subset\Delta$ . Then there exists a point $x \in
\Delta$ and an allocation function $\theta : \{0,1,\ldots,
d\}\rightarrow \{1,\ldots, q\}$ such that for each $\mu_i$ and all
$j=1,2,\ldots,q$,
\begin{equation}\label{eqn:fair-druga}
\sum_{\theta(\nu)=j} \mu_i({\rm cone}(x,F_\nu)) =
(1/q)\mu_i(\mathbb{R}^d)
\end{equation}
where ${\rm cone}(x,F_\nu)$ is the cone with the apex $x$
generated by the facet $F_i$.
\end{theo}

\medskip\noindent
{\bf Proof:} The proof is similar to the proof of
Theorem~\ref{thm:fair-prva} with the Theorem~\ref{thm:B-simplex}
used instead of Theorem~\ref{thm:randevu}. The details are left to
the reader. \hfill $\square$

\medskip
The proof of Theorem~\ref{thm:fair-druga} is somewhat simpler than
the proof of Theorem~\ref{thm:fair-prva} since the associated
configuration space is simpler. Nevertheless,
Theorem~\ref{thm:fair-prva} can be deduced from
Theorem~\ref{thm:fair-druga} by a limit and compactness argument.

\begin{prop}\label{prop:implication}
Theorem~\ref{thm:fair-druga} implies Theorem~\ref{thm:fair-prva}.
\end{prop}

\medskip\noindent
{\bf Proof:} (outline) Let us suppose that the supports of all
measures $\mu_i$ are contained in the ball $B(0,r)$, i.e.\ that
$\mu_i(\mathbb{R}^d\setminus B(0,r))=0$ for each $i=1,\ldots, n$.
Enlarge the simplex $\Delta$ by a homothety so that its diameter
$R$ is much larger than $r$, $R \gg r$. By
Theorem~\ref{thm:fair-druga} there exists a pair $(x,f)\in
\mathbb{B}(\Delta, S)$ which describes a fair dissection in the
sense of the equation (\ref{eqn:fair-druga}). There exists a
constant $C>1$ such that if $x\notin B(0, Cr)$ then for at least
one $\nu\in\{0,\ldots, d\}$ the cone ${\rm cone}(x, F_\nu)$ does
not intersect $B(0,r)$. It follows that if $(x,f)$ is a fair
division then $x\in B(0,Cr)$. Since $R \gg Cr$ we observe that the
cones ${\rm cone}(x, F_\nu)$ and $x + {\rm cone}(F_\nu)$ are very
close so by a limit and compactness argument the condition
(\ref{eqn:fair-prva}) can be also satisfied. \hfill $\square$

\section{Proofs of the polyhedral curtain theorem}

The `polyhedral curtain theorem' (Theorem~\ref{thm:prva}) is a
special case of Theorem~\ref{thm:fair-prva} for $q=2$, i.e. if
there are only two parties involved in the fair division. Indeed,
a polyhedral curtain is nothing but the common boundary of two
polyhedral sets,
\[
\bigcup_{\theta (\nu)=1}(x + {\rm cone}(F_\nu)) \quad \mbox{\rm
and} \quad \bigcup_{\theta (\nu)=2}(x + {\rm cone}(F_\nu)) .
\]
The reader interested only in the polyhedral curtain theorem may
use the Borsuk-Ulam theorem (along the lines of the proof of the
two dimensional case in Section~\ref{sec:two-dim-case}) together
with the $S=\{0,1\}$ case of Theorem~\ref{thm:randevu} which
claims that $$\mathbb{A}(R_\Delta, \mathcal{F}_\Delta, [2]) \cong
\partial(\diamondsuit^d)\cong S^{d-1}.$$
A further simplification can be achieved if one uses
Theorem~\ref{thm:B-simplex} together with the Borsuk-Ulam theorem,
and the limiting/compactness argument used in the proof of
Proposition~\ref{prop:implication}.

\section{Applications}

The `polyhedral curtain theorem' (Theorem~\ref{thm:prva}) is a
more combinatorial version od the ham sandwich theorem since it
involves alternatives.
\begin{figure}[hbt]
\centering
\includegraphics[scale=0.40]{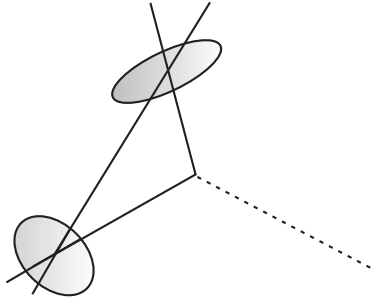}
\caption{Ham sandwich and polyhedral curtain equipartitions.}
\label{fig:comparison-1}
\end{figure}
For example in Figure~\ref{fig:comparison-1} we see an
equipartition of two measurable sets by a line and by one of three
planar `curtains' determined and prescribed in advance by the
chosen triangle $\Delta$.

It can be expected that some standard applications of the ham
sandwich theorem and its generalizations admit a modification
involving the polyhedral curtain theorem and its extensions. In
this section we offer only one example leaving a more complete
discussion for subsequent versions of this paper.

\subsection{Fair divisions by polynomial splines}

One of the standard and useful extensions (consequences) of the
ham sandwich theorem is the Stone-Tukey `polynomial ham sandwich
theorem', \cite{Stone-Tukey}. The key idea is to use some version
of the monomial Veronese embedding $V : \mathbb{R}^n \rightarrow
\mathbb{R}^N$ and read off the consequences in $\mathbb{R}^n$ of
the ham sandwich theorem applied in $\mathbb{R}^N$. Here we apply
a similar strategy to obtain consequences of the polyhedral
curtain theorem.

\medskip

For illustration we begin with an example.  Let $W : \mathbb{R}^2
\rightarrow \mathbb{R}^3$ be the embedding defined by $W(x,y) =
(x,y, x^2 + y^2)$ and let $\pi : \mathbb{R}^3\rightarrow
\mathbb{R}^2, \, \pi(x,y,z) = (x,y)$ be the projection. Suppose
that $X_1, X_2, X_3$ are measurable sets in $\mathbb{R}^2$ and let
$Y_i = W(X_i)$ be their images in the paraboloid $\Gamma\subset
\mathbb{R}^3$, described by the equation $z=x^2+y^2$.

An application of the polyhedral curtain theorem in $\mathbb{R}^3$
on the sets $\{Y_i\}_{i=1}^3$ yields a piecewise circular curve (a
circular spline) which cuts each of the sets $X_i$ in two parts of
equal measure. The number of nodes of the spline is controlled by
the number of faces in the polyhedral curtain.

\begin{figure}[hbt]
\centering
\includegraphics[scale=0.50]{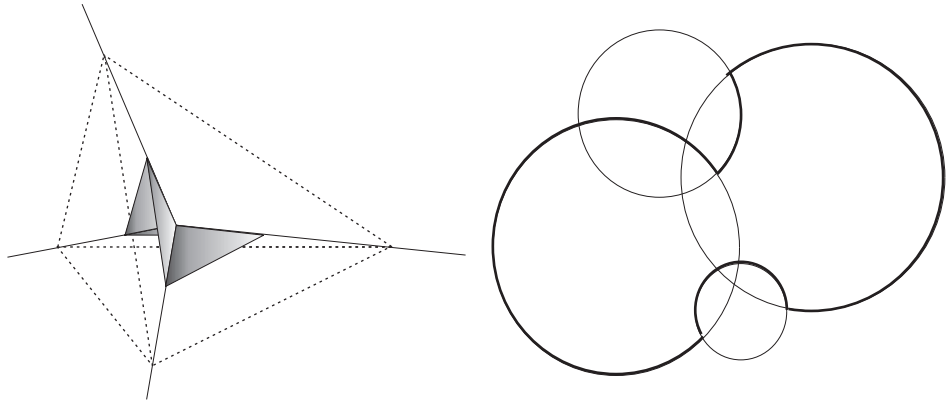}
\caption{A circular spline associated to a polyhedral curtain.}
\label{fig:curtain-6}
\end{figure}

Actually much more precise information about the spline can be
deduced from the simplex $\Delta\subset \mathbb{R}^3$. Indeed, the
slopes of the polyhedral faces of the curtain can be to some
extent prescribed in advance by the shape of the generating
tetrahedron $\Delta\subset \mathbb{R}^3$. Moreover, the
$\pi$-images $\pi(L_i\cap \Gamma)$ of the (non-empty)
intersections of two parallel planes $L_1$ and $L_2$ with the
paraboloid $\Gamma$ are two concentric circles in $\mathbb{R}^2$.

From here we deduce that the splines can be chosen to be
`concentric' to one of seven circular splines which are prescribed
in advance by the choice of the tetrahedron $\Delta$. These
splines are referred to as $\Delta$-generated.

\begin{prop}\label{prop:circular-splines}
Suppose that $X_1, X_2, X_3 \subset \mathbb{R}^2$ are three
measurable sets in the plane. Then for each simplex $\Delta\subset
\mathbb{R}^3$ there exists a $\Delta$-generated circular spline
which cuts in half each of the measurable sets $X_i$.
\end{prop}


\end{document}